\documentclass[reqno,draft,intlimits]{amsart}
\usepackage{amsfonts,amsmath,amsxtra,amsthm,amssymb,bezier,headerfo}

\theoremstyle{plain} %% This is the default

\newtheorem*{thm*}{Theorem}
\newtheorem*{cor*}{Corollary}
\newtheorem*{lem*}{Lemma}
\newtheorem*{prpr*}{Property}
\newtheorem*{prop*}{Proposition}
\newtheorem*{cond*}{Condition}
\newtheorem*{stat*}{Statement}

\newtheorem{thm}{Theorem}[section]

\newtheorem{lem}{Lemma}[section]

\theoremstyle{definition}

\newtheorem*{defn*}{Definition}
\newtheorem*{rem*}{Remark}
\newtheorem*{exm*}{Example}
\newtheorem*{exms*}{Examples}

\newtheorem{rem}{Remark}[section]

\allowdisplaybreaks

\newcommand{\R}{{\if mm {\rm I}\mkern -3mu{\rm R}\else \leavevmode
     \hbox{I}\kern -.17em \hbox{R} \fi}}
\renewcommand{\R}{{\mathbb R}}
\newcommand{\C}{{\if mm {{\rm C}\mkern -15mu{\phantom{\rm t}\vrule}}
    \mkern +10mu \else \leavemode \hbox{I}\kern -.17em \hbox{C} \fi}}
\renewcommand{\C}{{\mathbb C}}
\newcommand{\N}{{\if mm {\rm I}\mkern -3mu{\rm N}\else \leavevmode
\hbox{I}\kern -.17em \hbox{N} \fi}}
\renewcommand{\N}{{\mathbb N}}
\newcommand{\Z}{{\if mm {\sf Z}\mkern -8mu{\sf Z}\else \leavevmode
\hbox{Z}\kern -.17em \hbox{Z} \fi}}
\newcommand{\bB}{{\bf B}}

\newcommand{\bH}{{\bf H}}

\newcommand{\bR}{{\bf R}}

\newcommand{\cD}{{\mathcal D}}

\newcommand{\cH}{{\mathcal H}}

\newcommand{\ri}{{\rm i}}
\newcommand{\sD}{{\sf D}}

\newcommand{\sP}{{\sf P}}

\newcommand{\tA}{\widetilde{A}}
\newcommand{\ta}{\widetilde{a}}
\newcommand{\tb}{\widetilde{b}}

\newcommand{\D}{\displaystyle}

\newcommand{\Int}{\displaystyle\int\limits}
\newcommand{\Sup}{\mathop{\rm sup}}

\newcommand{\dist}{\mathop{\rm dist}}
\newcommand{\Img}{\mathop{\rm Im}}
\newcommand{\Real}{\mathop{\rm Re}}

\newcommand{\Var}{{\rm Var}}
\numberwithin{equation}{section}

\title{\sc Factorization Theorem for the Transfer Function
of a $2\!\times\!2$ Operator Matrix with Unbounded Couplings$^*$}

\author{V.~Hardt, R. Mennicken\\
University of Regensburg\\[2mm]
A. K. Motovilov\\
Joint Institute  for Nuclear Research, Dubna}

\address{V. Hardt, %R. Mennicken,
Department of Mathematics, University of Regensburg,
D-93040 Regensburg, Germany}

\email{volker.hardt@mathematik.uni-regensburg.de}

\address{R. Mennicken, Department of Mathematics, University of Regensburg,
D-93040 Regensburg, Germany}

\email{reinhard.mennicken@mathematik.uni-regensburg.de}

\address{A.K. Motovilov, Laboratory of Theoretical Physics,
Joint Institute  for Nuclear Research, Dubna, 141980, Russia}

\email{motovilv@thsun1.jinr.ru}

\thanks{Support of this work by the Deutsche Forschungsgemeinschaft,
the Heisenberg--Landau Program and the Russian Foundation for
Basic Research is gratefully acknowledged}  
\thanks{$^*$Contribution to Proceedings of the 1999 
Crimean Autumn Mathematical School}


\begin{document}
\oddpageheader{}{{\small Factorization Theorem for the Transfer
Function}}{{\small \thepage}}
\evenpageheader{{\small \thepage}}{{\small V. Hardt, R.
Mennicken, A. K. Motovilov}}{}

\begin{abstract}
We consider the analytic continuation of the transfer function
associated with a $2\times2$ operator matrix having unbounded
couplings into unphysical sheets of its Riemann surface.  We
construct a family of non-selfadjoint operators which factorize
the transfer function and reproduce certain parts of its
spectrum including the nonreal (resonance) spectrum situated in
the unphysical sheets neighboring the physical sheet.
\end{abstract}

\date{December 28, 1999}

\maketitle

%%%%%%%%%%%%%%%%%%%%%%%%%%%%%%%%%%%%%%%%%%%%%%%%%%%%%%%%%%%%%%%%%%%%%%
\section*{Introduction}\label{Intro}
%%%%%%%%%%%%%%%%%%%%%%%%%%%%%%%%%%%%%%%%%%%%%%%%%%%%%%%%%%%%%%%%%%%%%%
In this work we consider $2\times 2$ operator matrices
\begin{equation}\label{twochannel}
\bH_0\ =\ 
\left(\begin{array}{ccc}
A_0    & & T_{01} \\
T_{10} & &  A_1     \end{array}\right)
\end{equation}
acting in the orthogonal sum  $\cH=\cH_0\oplus\cH_1$ of separable
Hilbert spaces $\cH_0$ and $\cH_1$.  The entry
$A_0:\cH_0\rightarrow\cH_0$ is assumed to be an unbounded selfadjoint
operator with the domain  ${\mathcal D}(A_0)$.  We suppose that $A_0$
is semibounded from below, i.\,e., $A_0\geq\alpha_0$ for some
$\alpha_0\in\R $ and without loss of generality let $\alpha_0>0$. The
entry $A_1$ is assumed to be a bounded selfadjoint operator in $\cH_1$.
In contrast to~\cite{MennMotTMF,MennMotMathNachr},
in the present paper we
consider unbounded coupling operators $T_{ij}:\cH_j\rightarrow\cH_i$,
$i,j=0,1$, $i{\neq}j$.
Regarding these operators the following conditions
are supposed to be fulfilled:
\begin{equation}\label{IniCond}
T_{01}^*=T_{10}\quad\mbox{and}\quad
{\mathcal D}(T_{10})\supset{\mathcal D}(A_0^{1/2}).
\end{equation}
These assumptions are similar to those used in the works by {\sc
V.\,M.\,Adamyan, H.\,Lan\-ger, R.\,Mennicken} and {\sc
J.\,Sa\-u\-rer} \cite{AdLMSr}  and by {\sc R.\,Men\-ni\-cken}
and {\sc A.\,A.\,Shka\-likov} \cite{MenShk}.  Under the
conditions (\ref{IniCond}) the matrix~(\ref{twochannel}) is a
symmetric closable operator in $\cH$ on the domain ${\mathcal
D}(A_0)\oplus\cH_1$ and its closure $\bH=\overline{\bH}_0$ is a
selfadjoint operator (see \cite{AdLMSr,MenShk}); for an
explicit description of $\bH$ see \cite{AdLMSr,HMM,MenShk}.
Note that in applications arising from physical problems (see
e.\,g., Refs.~\cite{Goedbloed,Lifschitz,MotRem} and references
cited therein) one typically deals just with the case where
$\bH_0$ is a selfadjoint operator in a Hilbert space or a
symmetric operator admitting a selfadjoint closure.

The second condition in (\ref{IniCond}) yields that the product
$B_{10}:=T_{10}A_{0}^{-1/2}$ is a bounded linear operator.
It follows that $A_0^{-1/2}T_{01}$
has a bounded extension $B_{01}$ to the whole space $\cH_1$. This
extension coincides with $B_{10}^*$. In
addition, the hypothesis (\ref{IniCond}) implies that the
operator $A_1-z+T_{10}(z-A_0)^{-1}T_{01}$ for $z\in \varrho
(A_0)$ is densely defined and has a bounded extension onto the
whole space $\cH_1$ which is given by
\begin{equation}
\label{transferfunction}
M_1(z):=\tA _1-z+V_1(z)
\end{equation}
where $\tA_1:=A_1-B_{10}B_{01}$ and
$V_1(z):=zB_{10}(z-A_0)^{-1}B_{01}$.  We call $M_1(z)$ the
transfer function associated with the operator matrix $\bH$.  It
is obvious that this function, considered in the resolvent set
$\varrho (A_0)$ of the operator $A_0$, represents a holomorphic
operator-valued function.  (In the present work we use the
standard definition of holomorphy of an operator-valued
function with respect to the operator norm topology, see,
e.\,g., \cite{AdLMSr}).  It is worth noting that the holomorphic
operator-valued function $-M_1$ belongs to the class of
operator-valued Herglotz functions (see,
e.\,g., \cite{AkhiGlaz,GKMTs,KL1,Naboko}).

In the present paper like in \cite{MennMotTMF,MennMotMathNachr}
we study the transfer function $M_1(z)$ under the assumption
that it admits analytic continuation through the absolutely
continuous spectrum $\sigma_{\rm ac}(A_0)$ of the entry $A_0$.
We are especially interested in the case where the spectrum of
$\tA_1$ is partly or totally embedded into the absolutely
continuous spectrum of $A_0$. Notice that, since the resolvent
of the operator~$\bH$ can be expressed explicitly in terms of
$\bigl[M_1(z)\bigr]^{-1}$ (see, e.\,g.,
\cite{AdLMSr,MenShk,MennMotMathNachr}), in studying the spectral
properties of the transfer function one studies at the same time
the spectral properties of the operator matrix ${\bf H}$.

Section \ref{Transfer_function} includes a description of the 
conditions making analytic continuation of $M_1(z)$ through the 
spectrum $\sigma_{\rm ac}(A_0)$ possible.  Further, a 
representation of this continuation is given (see 
(\ref{Mcmpl})). In Section~\ref{SmainEq} we introduce the basic 
nonlinear equation (\ref{MainEq}) giving a rigorous sense to the 
formal operator equation $M_1(H_1)=0$. We explicitly show that 
eigenvalues and accompanying eigenvectors of a solution $H_1$ of 
the equation (\ref{MainEq}) are eigenvalues and eigenvectors of 
the analytically continued transfer function $M_1$. The 
solvability of (\ref{MainEq}) is proved under smallness 
conditions concerning the operator $B_{10}$, see (\ref{Best}).  
In Section~\ref{SecFactor} we first prove a factorization 
theorem (Theorem \ref{factorization}) for the analytically 
continued transfer function. This theorem implies that there 
exists certain domains in $\C$ lying partly on the unphysical 
sheet(s) where the spectrum of the analytically continued 
transfer function is represented by the spectrum of the 
corresponding solutions of the basic equation (\ref{MainEq}). 
Further in Section~\ref{SecFactor} we describe some relations 
between different solutions of (\ref{MainEq}) and some relations 
between their spectra.  Finally, in Section \ref{Example} we 
present a simple example.

A detailed exposition of the material presented here including proofs
in the case of essentially more general spectral situations
will be given in the extended paper \cite{HMM}.
%\newpage
%%%%%%%%%%%%%%%%%%%%%%%%%%%%%%%%%%%%%%%%%%%%%%%%%%%%%%%%%%%%%%%%%%%%%%
\section{Analytic continuation of the transfer function $M_1$}
\label{Transfer_function}
%%%%%%%%%%%%%%%%%%%%%%%%%%%%%%%%%%%%%%%%%%%%%%%%%%%%%%%%%%%%%%%%%%%%%%
For the sake of simplicity we assume in this paper that the
spectrum $\sigma (A_0)$ of the entry $A_0$ is absolutely
continuous consisting of the interval
$\Delta_0:=[\alpha_0,+\infty)$ with $\alpha_0>0$ while the
spectrum of $\tA _1$ is totally embedded into the interval
$\Delta _0$, i.\,e., $\sigma (\tA _1)\subset \Delta _0$.

Let  $E_0$ be the spectral measure for the entry $A_0$,
$A_0=\int_{\Delta_0}\lambda\,dE_0(\lambda)$.  Then the function
$V_1(z)$ can be written
$$
    V_1(z)\ =\ \int_{\Delta_0}dK_B(\mu)\frac{z}{z-\mu}
$$
with
$$
  K_B(\mu)\ :=\ B_{10}E^0(\mu)B_{01}
$$
where $E^0(\mu)$ stands for the spectral function of $A_0$,
$E^0(\mu)=E_0\bigl([\alpha_0,\mu)\bigr)$.
Thus, it is convenient to introduce the quantities
\begin{equation}\label{Varteta}
{\rm Var}_\theta(B)\ :=\ \Sup_{\left\{\delta_k,\,\mu_k\in\delta_k\right\}}
\sum_k (1+|\mu_k|)^{-\theta}\|B_{10}E_0(\delta_k)B_{01}\|,
\end{equation}
where $\theta$ is some real number and $\left\{\delta_k\right\}$
stands for a finite or countable complete system of Borel
subsets of $\sigma(A_0)=\Delta_0$ such that
$\delta_k\cap\delta_l=\emptyset$, if $k\neq l$, and
$\mathop{\bigcup}_k\delta_k=\Delta_0$. The points $\mu_k$ are
arbitrarily chosen points of $\delta_k$. The number $\mathop{\rm
Var}_\theta(B)$ is called weighted variation of the
operators $B_{ij}$ with respect to the spectral measure $E_0$.

Notice that in contrast to \cite{MennMotTMF,MennMotMathNachr},
where the variation~(\ref{Varteta}) was considered in case of
$\theta=0$, we now will mainly consider $\theta=1$. Of course,
introducing the variation ${\rm Var}_\theta(B)$ for $\theta\ne0$
only makes sense when the entry $A_0$ is an unbounded operator.

We suppose that the function $K_B(\mu)$ is differentiable in
$\mu\in\Delta_0$ in the operator norm topology.  The derivative
$K'_B(\mu)$ is non-negative, $ K'_B(\mu)\geq 0$,
since $K_B(\mu)$ is a non-decreasing function.
Obviously,
$$
   {\rm Var}_\theta(B)
   =\int_{\Delta_0}d\mu\,(1+|\mu|)^{-\theta}\|K'_B(\mu)\|.
$$

Further, we suppose that the function $K'_B(\mu)$ is continuous
within the interval $\Delta_0$ and, moreover, that it admits
analytic continuation from this interval to a simply
connected domain situated, say, in $\C^+$.
Let this domain be called $D^+$. We assume that
the boundary of the domain $D^+$ includes the entire spectral
interval $\Delta_0$. Since
$K'_B(\mu)$ represents a selfadjoint operator for
$\mu\in\Delta_0$ and $\Delta_0\subset\R$, the function
$K'_B(\mu)$ admits an analytic continuation from
$\Delta_0$ into the domain $D^-$, symmetric to $D^+$ with
respect to the real axis, $D^-=\{z:\,\overline{z}\in D^+\}$.
For the continuation into $D^-$ we will use the same notation
$K'_B(\mu)$. The selfadjointness of $K'_B(\mu)$ for
$\mu\in\Delta_{0}$ implies $[K'_B(\mu)]^*=K'_B(\bar{\mu}),$
$\mu\in D^\pm\,.$ Also, we shall suppose that the
$K'_B(\mu)$ satisfies the following condition at the
end point $\alpha_0$ of the spectral
interval $\Delta_0$:
$$
\|K'_B(\mu)\|\leq C|\mu-\alpha_0|^\gamma, \quad\mu\in D^\pm,
$$
with some $C>0$ and $\gamma\in (-1,0]$.

Let $\Gamma_{l}$ $(l=\pm 1)$ be a rectifiable Jordan curve in
$D^{l}$ resulting from continuous deformation of the interval
$\Delta_0$, the finite end point of this interval being fixed.
As mentioned above, in the following we deal with the variation
${\rm Var}_1(B)$. We extend the definition of this variation
also to the curve $\Gamma_l$ by introducing the modified
variation
\begin{equation}\label{NBnorm}
{\rm Var}_1(B,\Gamma_l)\,:=\, 
\displaystyle\int_{\Gamma_l}|d\mu|\,(1+|\mu|)^{-1}\|K'_B(\mu)\|
\end{equation}
where $|d\mu|$ denotes the Lebesgue measure on $\Gamma_l$.  We suppose that
the operators $B_{ij}$ are such that there exists a contour
(exist contours) $\Gamma _l$ on which the value ${\rm Var}_1(B,\Gamma_l)$
is finite, i.e., ${\rm Var}_1(B,\Gamma_l)<\infty$.
The contours $\Gamma_l$
satisfying the condition~${\rm Var}_1(B,\Gamma_l)<\infty$ are
said to be {$K_B$-bounded} contours. Surely, in the case of
unbounded $A_0$ the condition of boundedness of ${\rm
Var}_1(B,\Gamma_l)$ is much weaker than the condition of
boundedness of ${\rm Var}_0(B,\Gamma_l)$ used in
\cite{MennMotTMF,MennMotMathNachr}).
\begin{lem}\label{M1-Continuation}
The analytic continuation of the transfer function $M_1(z)$,
$z\in\C\setminus{\Delta_0}$,
through the spectral interval $\Delta _0$ into the subdomain
$D(\Gamma_l)\subset D^l$ $(l=\pm1)$ bounded by the set
$\Delta_0$ and
a $K_B$-bounded contour $\Gamma_l$ is given by
\begin{equation}\label{Mcmpl}
 M_1(z,\Gamma_l)\ :=\ \tA_1-z+V_1(z,\Gamma_l)
\end{equation}
where
\begin{equation}
\label{MGamma}
V_1(z,\Gamma_l)\ :=\ 
\int_{\Gamma_l} d\mu\,K'_B(\mu)\,\frac{z}{z-\mu}.
\end{equation}
For $z\in D^{l}\cap D(\Gamma_l)$ the function
$M_1(z,\Gamma_l)$ may be written as
\begin{equation}\label{M1Gresidue}
  M_1(z,\Gamma_l)\ =\ M_1(z)+2\pi\ri\,l z K'_B(z).
\end{equation}
\end{lem}
\begin{proof}
Obviously, the function~(\ref{MGamma}) is
well defined due to the $K_B$--boundedness of the contour
$\Gamma_l$ and
since for all $z\in\C\setminus \Gamma_l$
there exist a  $c(z)>0$ such that the estimate
$\bigl |(z-\mu)^{-1}\bigr |<c(z)\,\bigl (1+|\mu|\bigr )^{-1}$
$(\mu\in \Gamma _l\bigr )$ holds.
Thus, the proof of
this lemma is reduced to the observation that the
function $M_1(z,\Gamma_l)$ is holomorphic for
$z\in\C\setminus \Gamma_l$ and coincides with
$M_1(z)$ for
$z\in\C\setminus \overline{D(\Gamma_l)}$.
The equation~(\ref{M1Gresidue}) is obtained from~(\ref{MGamma}) using the
Residue Theorem.
\end{proof}

The formula~(\ref{M1Gresidue}) shows that in general the
transfer function $M_1(z)$ has a Riemann
surface with at least two sheets.  The sheet of the complex plane
where the transfer
function $M_1(z)$ together with the resolvent
$\bR(z)=(\bH -z)^{-1}$ is initially considered  is said to be the {physical
sheet}.  The remaining sheets of the Riemann surface of $M_1(z)$
are said to be unphysical sheets (see, e.\,g.,
\cite{ReedSimonIII}). In the present work we deal with the
unphysical sheets neighboring the physical one, i.\,e., with the
sheets connected through the interval $\Delta_0$ immediately to
the physical sheet.
\begin{rem}\label{r2.2}
For $z\in\C\setminus\Gamma_l$, the equation~{\rm (\ref{MGamma})} defines
values of the function $V_1(\cdot,\Gamma_l)$ in the space of
bounded operators in $\cH_1$. The inverse
transfer function $\bigl [M_1(z)\bigr ]^{-1}$ coincides with the
right lower block component $\bR_{11}(z)$ of the resolvent
$\bR(z)=(\bH-z)^{-1}$ and, thus, it is holomorphic in
$\C\setminus\sigma(\bH)\supset\C\setminus\R$. Since
$M_1(z,\Gamma_l)$ coincides with $M_1(z)$ for all
$z\in\C\setminus \overline{D(\Gamma_l)}$, one concludes that
$[M_1(z,\Gamma_l)]^{-1}$ exists as a bounded operator and is
holomorphic in $z$ at least for
$z\in\C\setminus\bigl(\sigma(\bH)\cup\overline{D(\Gamma_l)}\bigr)$.
\end{rem}
%%%%%%%%%%%%%%%%%%%%%%%%%%%%%%%%%%%%%%%%%%%%%%%%%%%%%%%%%%%%%%
\section{The basic equation}
\label{SmainEq}
%%%%%%%%%%%%%%%%%%%%%%%%%%%%%%%%%%%%%%%%%%%%%%%%%%%%%%%%%%%%%%
Let $\Gamma$ be a $K_B$-bounded contour.  If $Y$ stands for an
arbitrary bounded operator in $\cH_1$ such that the spectrum of
$Y$ is separated from the set $\Gamma$ then, following
to~\cite{MennMotMathNachr,MotSPbWorkshop,MotRem},
one can define the  operator
\begin{equation}\label{V1Y}
V_1(Y,\Gamma)\ :=\ \int_{\Gamma}
d\mu\,K'_B(\mu)\,Y(Y-\mu)^{-1}.
\end{equation}
This operator is bounded, $V_1(Y,\Gamma)\in\bB(\cH_1,\cH_1)$,
and its norm admits the estimate
\begin{equation}
\label{V1Yest}
\|V_1(Y,\Gamma)\|\ \leq\ {\rm Var}_1(B,\Gamma)\ 
\|Y\|\ \Sup_{\mu\in\Gamma}(1+|\mu|)\|(Y-\mu)^{-1}\|.
\end{equation}
In what follows we consider the equation
(cf. \cite{MennMotMathNachr,MotSPbWorkshop,MotRem})
\begin{equation}\label{MainEq}
Y\ =\ \tA_1+V_1(Y,\Gamma).
\end{equation}
This equation possesses the following
characteristic property: If an operator $H_1$ is a
solution of~(\ref{MainEq}) and $u_1$ is an eigenvector of $H_1$,
$H_1 u_1=zu_1$, then
\begin{eqnarray*}
zu_1&=&\tA_1 u_1+V_1(H_1,\Gamma )u_1\ =\ 
\tA_1u_1+\int_{\Gamma}d\mu\,K'_B(\mu )H_1(H_1-\mu)^{-1}u_1\\
&=&\tA_1u_1+\int_{\Gamma}d\mu\,K'_B(\mu )\frac z{z-\mu}u_1
\ =\ \tA_1 u_1+V_1(z,\Gamma)u_1.
\end{eqnarray*}
This means that any eigenvalue $z$ of such an operator $H_1$
is automatically an eigenvalue for the analytically continued transfer
function $M_1(z,\Gamma )$ and $u_1$ is a corresponding eigenvector.  Thus,
having found the solution(s) of the equation (\ref{MainEq}) one
obtains an effective means of studying the spectral properties
of the transfer function $M_1(z,\Gamma)$, referring to
well-known facts of operator theory \cite{GK,Kato}.  It is
convenient to rewrite the equation~(\ref{MainEq}) in the form
\begin{equation}\label{MainEqC}
X\ =\ V_1(\tA_1+X,\Gamma)
\end{equation}
where $X:=Y-\tA_1$.

Let the spectrum of the operator
$\tA_1$ be separated from $\Gamma$, i.\,e.,
\begin{equation}\label{dist-tA}
d_0(\Gamma)\ :=\ \dist\{\sigma(\tA_1),\Gamma\}>0\,.
\end{equation}
Then, since $\tA_1$ is selfadjoint and bounded, it is obvious
that the following quantity
\begin{equation}\label{VarBtAGdef}
{\rm Var}_{\tA_1}(B,\Gamma)\ := \
\int_{\Gamma}|d\mu|\,\frac{\|K'_B(\mu)\|}{\dist\{\mu,\sigma(\tA_1)\}}
\end{equation}
is finite,
\begin{equation}\label{VarBtAGfin}
{\rm Var}_{\tA_1}(B,\Gamma)\ \leq\ \Var_1(B,\Gamma)\,\mathop{\rm sup}
 \limits_{\mu\in\Gamma}
(1+|\mu|)\bigl [\dist\{\mu,\sigma(\tA_1)\}\bigr ]^{-1}\ <\ \infty\,.
\end{equation}
It is more convenient to make the subsequent estimations in terms
of the variation ${\rm Var}_{\tA_1}(B,\Gamma)$ rather then in terms
of the variation ${\rm Var}_1(B,\Gamma)$.
\begin{thm}\label{Solvability}
Let $\tA_1$ be a bounded operator, the
contour $\Gamma$ be $K_B$-bounded and
\begin{equation}\label{Best}
\Var_{\tA_1}(B,\Gamma)<1\,, \qquad
\Var_{\tA_1}(B,\Gamma)\|\tA_1\|<
\displaystyle\frac{1}{4}\,d_0(\Gamma)\,[1-\Var_{\tA_1}(B,\Gamma)]^2\,.
\end{equation}
Let
\begin{equation}\label{rmin}
\begin{array}{rcl}
\displaystyle r_{\rm min}(\Gamma)&:=&\displaystyle \frac{1}{2}\,d_0(\Gamma)\,
[1-\Var_{\tA_1}(B,\Gamma)]\\[3mm]
&&\displaystyle -\sqrt{\frac{1}{4}\,d^2_0(\Gamma)\,[1-\Var_{\tA_1}(B,\Gamma)]^2
-d_0(\Gamma)\,\Var_{\tA_1}(B,\Gamma)\,\|\tA_1\|}\quad
\end{array}
\end{equation}
and
\begin{equation}\label{rmax}
r_{\rm max}(\Gamma) \ :=\  d_0(\Gamma)-\sqrt{\Var_{\tA_1}(B,\Gamma)\,
d_0(\Gamma)\,[d_0(\Gamma)+\|\tA_1\|]}.
\end{equation}
Then the equation~{\rm(\ref{MainEqC})} is uniquely solvable in any
closed ball
$${\mathcal S}_1(r)\ :=\ \bigl\{X\in \bB(\cH_1,\cH_1)\, :\,
\Vert X\Vert \leq r\bigr\}$$
where
\begin{equation}\label{Br}
  r_{\rm min}(\Gamma)\ \leq \ r\ < \ r_{\rm max}(\Gamma).
\end{equation}
The solution $X$ of the equation~{\rm(\ref{MainEqC})} is the same for any
$r$ satisfying~{\rm(\ref{Br})} and in fact it belongs to the smallest
ball ${\mathcal S}_1(r_{\rm min})$, $\|X\|\leq r_{\rm min}(\Gamma)$.
\end{thm}
\begin{proof}
One can prove this theorem making use of Banach's Fixed Point Theorem
(see \cite{HMM}).
\end{proof}
The following statement is a direct consequence of the conditions
(\ref{Best}).
\begin{rem}
\label{Half1mVar}
The values of $r_{\rm min}(\Gamma)$ and $r_{\rm max}(\Gamma)$
satisfy the estimates
$$
   r_{\rm min}(\Gamma)\ <\ \frac{1}{2}\,
   d_0(\Gamma)\,[1-\Var_{\tA_1}(B,\Gamma)]
   \ <\ r_{\rm max}(\Gamma)\,.
$$
\end{rem}
\begin{thm}\label{Hunique}
Let the conditions of {\rm Theorem~\ref{Solvability}} be fulfilled
for a $K_B$-bounded contour $\Gamma\subset D^l$ and let $X$ be
the solution of the equation{~\rm(\ref{MainEqC})}. Then
$X$ coincides with the analogous solution $\widetilde{X}$ for
any other $K_B$-bounded contour $\widetilde{\Gamma}\subset D^l$
satisfying the estimates
$$
\Var_{\tA_1}(B,\tilde \Gamma)<1\quad and\quad
\Var_{\tA_1}(B,\tilde \Gamma)\|\tA_1\|<
\displaystyle\frac{1}{4}\,\widetilde{d}_0
[1-\Var_{\tA_1}(B,\tilde \Gamma)]^2
$$
where
$0 < \widetilde{d}_0=\dist\{\sigma(\tA_1),
\sigma'(A_0)\cup\widetilde{\Gamma}\}\leq d_0(\Gamma)$.
Moreover, this solution satisfies the inequality
$\|X\|\leq r_0(B)$ where
$$
  r_0(B):=\inf \bigl\{r_{\rm min}(\Gamma_l)\, :\,
  \Var_{\tA_1}(B,\Gamma_l)<1\,,  \,\omega(B,\Gamma_l)>0\bigr\}
$$
with $r_{\rm min}(\Gamma_l)$ given by~{\rm(\ref{rmin})} and
$$
\omega(B,\Gamma_l)\ :=\ d_0(\Gamma _l)\,[1-\Var_{\tA_1}(B,\Gamma_l)]^2
-4\|\tA_1\|\Var_{\tA_1}(B,\Gamma_l).
$$
The value of $r_0(B)$ does not depend on $l$.
\end{thm}

So, for a given holomorphy domain $D^l$ $(l=\pm 1)$ the solutions $X$ and
$H_1$, $H_1=\tA_1+X,$ do not depend on the $K_B$-bounded contours
$\Gamma_l\subset D^l$ satisfying the conditions~(\ref{Best}).
But when the index $l$ changes,  $X$ and $H_1$ can also change.
For this reason we shall supply them in the following, when it
is necessary, with the index $l$ writing
$X^{(l)}$ and $H_1^{(l)}$, $H_1^{(l)}=\tA_1+X^{(l)}$.
Surely, the equations~(\ref{MainEq}) and~(\ref{MainEqC}) are nonlinear
equations and, outside the balls $\|X\|<r_{\rm max}(\Gamma_l)$,
they may, in principle, have other solutions, different from the
$X^{(l)}$ or $H_1^{(l)}$ the existence of which is guaranteed by
Theorem~\ref{Solvability}. In the following we only deal with
the solutions $X^{(l)}$ or $H_1^{(l)}$ for $l=\pm 1$.
%%%%%%%%%%%%%%%%%%%%%%%%%%%%%%%%%%%%%%%%%%%%%%%%%%%%%%%%%%%%%%%%%%%%
\section{Factorization theorem}\label{SecFactor}
%%%%%%%%%%%%%%%%%%%%%%%%%%%%%%%%%%%%%%%%%%%%%%%%%%%%%%%%%%%%%%%%%%%%
Now we prove a {\it factorization theorem} for the transfer
function $M_1(z,\Gamma_l)$. Note that this theorem recalls the
corresponding statements from \cite{MarkusMatsaev,VirozubMatsaev}.
\begin{thm}\label{factorization}
Let $\Gamma_l$ be a $K_B$-bounded contour satisfying the
conditions~{\rm(\ref{Best})}.
Suppose $X^{(l)}$ is the solution of the basic
equation~{\rm(\ref{MainEqC})}, $\|X^{(l)}\|\leq r_0(B)$, and
$H_1^{(l)}=\tA_1+X^{(l)}$. Then,
for $z\in\C\setminus \Gamma_l$, the
transfer function $M_1(z,\Gamma_l)$ admits the factorization
\begin{equation}\label{Mfactor}
    M_1(z,\Gamma_l)\ =\ W_1(z,\Gamma_l)\,(H_1^{(l)}-z)\,
\end{equation}
where $W_1(z,\Gamma_l)$ is a bounded operator in $\cH_1$,
\begin{equation}\label{Mtild}
\begin{array}{rcl}
\displaystyle W_1(z,\Gamma_l)&=&\displaystyle I_1-\int_{\Gamma_l}
d\mu\,K'_B(\mu)\,(H_1^{(l)}-\mu)^{-1}\\[7mm]
&&\displaystyle +z\,\int_{\Gamma_l}
d\mu\,K'_B(\mu)(z-\mu)^{-1}(H_1^{(l)}-\mu)^{-1}\,.
\end{array}
\end{equation}
If $\dist\{z,\sigma(\tA_1)\}\leq
d_0(\Gamma_l){[1-\Var_{\tA_1}(B,\Gamma_l)]/2}$, then the operator
$W_1(z,\Gamma_l)$ is boundedly invertible and
\begin{equation}\label{Mtest}
  \left\|[W_1(z,\Gamma_l)]^{-1}\right\|
  \ \leq \ \left(1-\frac{4\,\Var_{\tA_1}(B,\Gamma_l)\,
  \bigl[d_0(\Gamma_l)+\|\tA_1\|\bigr]}
  {d_0(\Gamma_l)\bigl[1+\Var_{\tA_1}(B,\Gamma_l)\bigr]^2}\right)^{-1}\ <
  \ \infty.
\end{equation}
\end{thm}
\begin{proof}
First we prove the formula~(\ref{Mfactor}).
Note that, according to~(\ref{V1Y}) and~(\ref{MainEqC}),
\begin{equation}\label{e4.4}
\tA_1\ =\ H_1^{(l)}-V_1(\tA_1+X^{(l)},\Gamma_l)
\ =\ H_1^{(l)}-\int_{\Gamma_l}
d\mu\,K'_B(\mu )H_1^{(l)}(H_1^{(l)}-\mu )^{-1}\,.
\end{equation}
Thus, in view of the representations~(\ref{Mcmpl}) and~(\ref{MGamma}),
the function $M_1(z,\Gamma_l)$ can be written as
\begin{eqnarray*}
M_1(z,\Gamma_l)&=&\tA_1-z+\int_{\Gamma_l}
d\mu\,K'_B(\mu )\frac z{z-\mu }\\
&=&H_1^{(l)}-z-\int_{\Gamma_l}
d\mu\,K'_B(\mu)(H_1^{(l)}-\mu )^{-1}\bigl(H_1^{(l)}-z\bigr)\\
&&+z\int_{\Gamma_l}
d\mu\,K'_B(\mu)\left[ \frac {1}{z-\mu}-(H_1^{(l)}-\mu )^{-1}\right]\\
&=&\bigl (H_1^{(l)}-z\bigr )-\int_{\Gamma_l}
d\mu\,K'_B(\mu)(H_1^{(l)}-\mu )^{-1}\bigl(H_1^{(l)}-z\bigr)\\
&&+z\int_{\Gamma_l}
d\mu\,K'_B(\mu) (z-\mu)^{-1}(H_1^{(l)}-\mu )^{-1} \bigl(H_1^{(l)}-z\bigr).
\end{eqnarray*}
which proves the equation (\ref{Mfactor}).
The boundeness of the operator $W_1(z,\Gamma_l)$
for $z\in\C\setminus \Gamma_l$
is obvious.

Further, we give a sketch of the proof that the factor $W_1(z,\Gamma_l)$ is a
boundedly invertible operator if the condition
$\dist\{z,\sigma(\tA_1)\}
\leq{d_0(\Gamma_l)\,[1-\Var_{\tA_1}(B,\Gamma_l)]/2}$ holds.
The formula
\begin{equation}\label{ResEstim}
\|(\tA_1+X^{(l)}-\mu)^{-1}\|
\ \leq \
\frac{1}{\dist\{\mu,\sigma(\tA_1)\}-\|X\|}\,,
\end{equation}
the definitions of $d_0(\Gamma_l)$ and
$r_{\min}(\Gamma _l)$ and Remark \ref{Half1mVar} imply that
\begin{equation}\label{Est00}
\left\|\ \int_{\Gamma_l}
d\mu\,K'_B(\mu)\,(H_1^{(l)}-\mu)^{-1}\right\| \ \leq \
\frac{2\,\Var_{\tA_1}(B,\Gamma_l)}{1+\Var_{\tA_1}(B,\Gamma_l)}\,.
\end{equation}
Using again the inequality (\ref{ResEstim}) and Remark
\ref{Half1mVar} we find
$$
\left\|z\int_{\Gamma_l}
d\mu\,K'_B(\mu)\,(H_1^{(l)}-\mu)^{-1}(z-\mu)^{-1}\right\|
\ \leq\ |z|\,\frac{2\,\Var_{\tA_1}(B,\Gamma_l)}{1+\Var_{\tA_1}(B,\Gamma_l)}\,
\sup\limits_{\mu\in \Gamma_l}|z-\mu|^{-1}\,.
$$
The inequality
$\dist\{z,\sigma(\tA_1)\}\leq{d_0(\Gamma_l)[1-\Var_{\tA_1}(B,\Gamma_l)]/2}
$
yields
$$
|z|\ \leq\ \|\tA_1\|+\dist\bigl \{z,\sigma (\tA_1)\bigr \}\leq
\|\tA_1\|+\frac{1}{2}d_0(\Gamma_l)[1-\Var_{\tA_1}(B,\Gamma_l)]
$$
and one obtains for $\mu\in \Gamma_l$ that
$$
\sup\limits_{\mu\in \Gamma_l}
|z-\mu|^{-1}\ \leq\ \frac{2}{d_0(\Gamma_l)[1+\Var_{\tA_1}(B,\Gamma_l]}.
$$
Hence, for $\dist\{z,\sigma(\tA_1)\}\leq{d_0(\Gamma_l)
[1-\Var_{\tA_1}(B,\Gamma_l)]/2}$,
\begin{eqnarray*}
\|W_1(z,\Gamma_l)-I_1\| &\leq&
\frac{2\,\Var_{\tA_1}(B,\Gamma_l)}{1+\Var_{\tA_1}(B,\Gamma_l)} \\
&& +\frac{4\,\Var_{\tA_1}(B,\Gamma_l)\,
\left\{\|\tA_1\|+\frac{1}{2}\,d_0(\Gamma_l)
[1-\Var_{\tA_1}(B,\Gamma_l)]\right\}}
{d_0(\Gamma_l)\,[1+\Var_{\tA_1}(B,\Gamma_l)]^2} \\
&=& \frac{4\,\Var_{\tA_1}(B,\Gamma_l)\,\,[d_0(\Gamma_l)+\|\tA_1\|]}
{d_0(\Gamma_l)\,[1+\Var_{\tA_1}(B,\Gamma_l)]^2}\ <\ 1\,.
\end{eqnarray*}
The last inequality is a direct consequence
of the second assumption in (\ref{Best}). We conclude
that $W_1(z,\Gamma_l)$ is invertible and
that the inequality (\ref{Mtest}) holds.
\end{proof}
The following theorems can be proved in the same way as
Theorem~4.4  and Theorem~ 4.7 in \cite{MennMotMathNachr}.
\begin{thm}\label{SpHalfVic}
The spectrum $\sigma(H_1^{(l)})$ of the operator
$H_1^{(l)}=\tA_1+X^{(l)}$ belongs to the closed $r_0(B)$-neighbourhood
$$
{\mathcal O}_{r_0(B)}(\tA_1)\ :=\ \bigl\{z\in\C:\,\dist\{z,\sigma(\tA_1)\}
\leq r_0(B)\bigr \}
$$
of the spectrum of $\tA_1$.  If the contour $\Gamma_l\subset 
D^l$ satisfies~{\rm(\ref{Best})}, then the nonreal spectrum of 
$H_1^{(l)}$ belongs  to $D^l\cap{\mathcal O}_{r_0(B)}(\tA_1)$. 
Moreover, the spectrum $\sigma(H_1^{(l)})$ coincides with a 
subset of the spectrum of the transfer function 
$M_1(\cdot,\Gamma_l)$. More precisely, the spectrum of 
$M_1(\cdot,\Gamma_l)$ in the set
$$
{\mathcal O}(\tA_1,\Gamma_l)\ :=\ \left\{z\in\C\, : \,
\dist\{z,\sigma(\tA_1)\}\leq
d_0(\Gamma_l)\,[1-\Var_{\tA_1}(B,\Gamma_l)]/{2}\right\}
$$
equals the spectrum of $H_1^{(l)}$, i.\,e.,
\begin{equation}\label{e3.7new}
\sigma\bigl (M_1(\cdot,\Gamma_l)\bigr)\cap {\mathcal
O}(\tA_1,\Gamma_l)\,=\,\sigma(H_1^{(l)}).
\end{equation}
In fact such a statement
separately holds for the point and continuous spectra.
\end{thm}

In the following lemma we state a simple but useful relation between
$H_1^{(l)}$ and the adjoint operator of $H_1^{(-l)}$.
According to our convention
$\Gamma_{(-l)}\subset D^{(-l)}$ is the contour
which is conjugate to the contour $\Gamma_l$.
\begin{lem}\label{Adjoint}
Let $\Gamma_l\subset D^l$ be a $K_B$-bounded contour
for which the conditions of {\rm Theorem~\ref{Solvability}}
are fulfilled. Then for any
$z\in\C\setminus\Gamma_l$
the following equality holds true:
\begin{equation}\label{Hadj}
W_1(z,\Gamma_l)\,\left (H_1^{(l)}-z\right)\ =\ 
\left( H_1^{(-l)*}-z\right)\,
\left [W_1(\overline{z},\Gamma_{(-l)})\right ]^*\,.
\end{equation}
Further the spectrum of $H_1^{(-l)*}$
coincides with the spectrum of $H_1^{(l)}$.
\end{lem}
\begin{proof}
Let $z\in \C\setminus \Gamma_l$. By definition $\overline{z}\in
\C\setminus \Gamma_{(-l)}$ and
\begin{equation}\label{e3.8new}
M_1(z,\Gamma_l)^* \ = \ M_1(\overline{z},\Gamma_{(-l)}).
\end{equation}
Therefore, the relation (\ref{Hadj}) follows from 
the factorizations
$$
M_1(z,\Gamma_l)\,= \,W_1(z,\Gamma_l)\, (H_1^{(l)}-z)
$$
and
$$
M_1(\overline{z},\Gamma_{(-l)})\,=\,
W_1(\overline{z},\Gamma_{(-l)})\, (H_1^{(-l)}-\overline{z}).
$$
By the relation (\ref{e3.7new}) $z$ belongs to the spectrum of the operator 
$H_1^{(-l)*}$ if and only if 
$\overline{z}\in {\mathcal O}(\tA_1,\Gamma_{(-l)})$ and 
$0\in\sigma\bigl([M_1(\overline{z},\Gamma_{(-l)})]^*\bigr)$. {}From 
(\ref{e3.8new}) we conclude that 
$0\in\sigma\bigl([M_1(\overline{z},\Gamma_{(-l)})]^*\bigr)$ if and
only if $0\in \sigma \bigl (M_1(z,\Gamma _{l})\bigr )$. Again by
(\ref{e3.7new}) the coincidence of the spectra of $H_1^{(l)}$ and 
$H_1^{(-l)*}$ follows.
\end{proof}
Let
\begin{equation}\label{e4.8neu}
\Omega^{(l)}\ :=\ \displaystyle\int_{\Gamma_l}
d\mu\,\mu\,\,(H_1^{(-l)*}-\mu)^{-1}K'_B(\mu )\,(H_1^{(l)}-\mu)^{-1}\,
\end{equation}
where as previously $\Gamma_l$ denotes a $K_B$-bounded
contour satisfying the conditions~(\ref{Best}).  The operator
$\Omega^{(l)}$ does not depend on the choice
of such a $\Gamma_l$.
%%%%%%%%%%%%%%%%%%%%%%%%%%%%%%%%%%%%%%%%%%%%%%%%%%%%%%%%%%%%%%%%%%%
\begin{thm}\label{MHOmega}
The operators $\Omega^{(l)}$ $(l=\pm 1)$ possess the following properties
{\rm(}cf.
{\rm\cite{HMM,MarkusMatsaev,MennMotTMF,MennMotMathNachr,MenShk,VirozubMatsaev}):}
\begin{eqnarray}
\label{Omest}
\|\Omega^{(l)}\|<1, \qquad  \Omega^{(-l)}&=&\Omega^{(l)*},\\[2mm]
\label{MOmega}
-\frac{1}{2\pi\ri}\int_\gamma dz\,[M_1(z,\Gamma_l)]^{-1} &=&
(I_1+\Omega^{(l)})^{-1}\,,\\[2mm]
\label{HOmega}
-\frac{1}{2\pi\ri}\int_\gamma dz\,z\,[M_1(z,\Gamma_l)]^{-1} &=&
(I_1+\Omega^{(l)})^{-1}H_1^{(-l)*}=
H_1^{(l)}(I_1+\Omega^{(l)})^{-1}\,,\qquad
\end{eqnarray}
where $\gamma$ stands for an arbitrary rectifiable closed
contour going around the spectrum of $H_1^{(l)}$ inside the set ${\mathcal
O}(\tA_1,\Gamma_l)$ in the positive direction.  
The integration along $\gamma$ is
understood in the strong sense.
\end{thm}
\begin{proof}
The estimate in (\ref{Omest}) can be proved by using 
the relation (\ref{Hadj}) following 
the proof of the estimate (\ref{Mtest}).
This estimate yields that the sum
$I_1+\Omega^{(l)}$ is a boundedly invertible
operator in $\cH_1$.

To prove the formula~(\ref{MOmega}) we recall that  
due to the factorization theorem~\ref{factorization}
and the formula~(\ref{Hadj}) the following
factorization holds for 
for $z\in{\mathcal O}(\tA_1,\Gamma_l)\backslash \sigma (H_1^{(l)})$: 
\begin{equation}\label{e4.12z}
\begin{array}{rcl}
[M_1(z,\Gamma_l)]^{-1}& = &\left (H_1^{(l)}-z\right)^{-1}\,
[W_1(z,\Gamma_l)]^{-1}\\[3mm]
&=&[W_1(\overline{z},\Gamma_{(-l)})]^{*-1}\,
\left(H_1^{(-l)*}-z\right)^{-1}
\end{array}
\end{equation}
where $[W_1(z,\Gamma_l)]^{-1}$ and $[W_1(\overline{z},\Gamma_{(-l)})]^{*-1}$
are holomorphic functions with values in $\bB(\cH_1,\cH_1)$.
By the resolvent equation and the definition (\ref{Mtild}) 
the product $\Omega^{(l)}(H_1^{(l)}-z)^{-1}$
can be written as
\begin{equation}\label{e4.12zz}
\Omega^{(l)}(H_1^{(l)}-z)^{-1}\ =\ F_1(z)+F_2(z)
\end{equation}
where
\begin{equation}\label{F1int}
F_1(z)\ :=\ \int_{\Gamma_l}
d\mu\,\mu\,(H_1^{(-l)*}-\mu)^{-1}\,K'_B(\mu)\,(H_1^{(l)}-\mu)^{-1}
\, (\mu-z) ^{-1}
\end{equation}
and
\begin{equation}\label{e4.13z}
\begin{array}{rcl}
\displaystyle F_2(z)&:=&
\displaystyle\left (-\int_{\Gamma_l}d\mu
\frac {\mu}{\mu -z}\,(H_1^{(-l)*}-\mu)^{-1}\,K'_B(\mu)\right )\,
(H_1^{(l)}-z)^{-1}\\[3mm]
&=&
\displaystyle\left([W_1(\overline{z},\Gamma_{(-l)})]^*-I_1\right )
(H_1^{(l)}-z)^{-1}.
\end{array}
\end{equation}
Further, the formula (\ref{e4.12z})
yields that
$$
(I_1+\Omega^{(l)})\,[M_1(z,\Gamma_l)]^{-1}\ = \
F_1(z)\,[W_1(z,\Gamma_l)]^{-1}+(H_1^{(-l)*}-z)^{-1}.
$$
The function $F_1(z)$ is holomorphic inside the contour
$\gamma$, $\gamma\subset{\mathcal O}(\tA_1,\Gamma_l)$, since the
argument $\mu$ of the integrand in the formula~(\ref{F1int}) belongs to
$\Gamma_l$ and  thereby
\mbox{$|z-\mu|\geq [d_0(\Gamma_l)+\Var_{\tA_1}(B,\Gamma_l)]/2>0$.} Thus
the term $F_1(z)[W_1(z,\Gamma_l)]^{-1}$ does not contribute to
the integral
$$
-\D\frac{1}{2\pi\ri}\Int_\gamma dz (I_1+\Omega^{(l)})
[M_1(z,\Gamma_l)]^{-1}
$$
while the resolvent $(H_1^{(-l)*}-z)^{-1}$
gives the identity $I_1$ which proves the equation~(\ref{MOmega}). 

Regarding the equation~(\ref{HOmega}) we obtain
\begin{eqnarray*}
   \lefteqn{-\D\frac{1}{2\pi\ri}\Int_\gamma dz (I_1+\Omega^{(l)})
   \,z\, [M_1(z,\Gamma_l)]^{-1} =}\\
  &=&-\D\frac{1}{2\pi\ri}\Int_\gamma dz
  \,z\, F_1(z)\,[W_1(z,\Gamma_l)]^{-1}
  -\D\frac{1}{2\pi\ri}\Int_\gamma dz \,z\,(H_1^{(-l)*}-z)^{-1}\,.
\end{eqnarray*}
The first integral vanishes whereas the second integral equals 
$H_1^{(-l)*}$. The second equation of~(\ref{HOmega})
can be checked in the same way.
\end{proof}

Note that the formulae~(\ref{MOmega}) and~(\ref{HOmega}) allow, in
principle, to construct the operators $H_1^{(l)}$, $l=\pm 1$, and, thus, to
resolve the equation~(\ref{MainEqC}) by a contour integration of
the inverse of the transfer function $M_1(z,\Gamma _l)$.
\begin{rem}\label{HlHml}
The formula~{\rm(\ref{HOmega})} implies that
$$
H_1^{(l)*}\ =\ 
(I_1+\Omega^{(-l)})\,H_1^{(-l)}\,(I_1+\Omega^{(-l)})^{-1}\,.
$$
Therefore the spectrum of $H_1^{(-l)*}$
coincides with the spectrum of $H_1^{(l)}$.
\end{rem}
\begin{thm}
\label{MPOmega}
Let $\lambda$ be an isolated eigenvalue of the operator
$H_1^{(l)}$ and, consequently, of the operator $H_1^{(-l)*}$ and
of the transfer function $M_1(z,\Gamma_l)$ taken for a
$K_B$-bounded contour $\Gamma_l$ satisfying the
conditions~{\rm(\ref{Best})}. By $\sP_\lambda^{(l)}$ and
$\sP_\lambda^{(-l)*}$ we denote the eigenprojections of the
operators $H_1^{(l)}$ and $H_1^{(-l)*}$, respectively, and by
$P_\lambda^{(l)}$ the residue of $M_1(z,\Gamma_l)$ at
$z=\lambda$,
\begin{eqnarray}\label{Plambda}
\sP_\lambda^{(l)}& :=& 
-\D\frac{1}{2\pi\ri}\Int_\gamma dz\,\,(H_1^{(l)}-z)^{-1},
\\
\label{Plambdastar}
\sP_\lambda^{(-l)*}& :=& 
-\D\frac{1}{2\pi\ri}\Int_\gamma dz\,\,(H_1^{(-l)*}-z)^{-1}
\end{eqnarray}
and
\begin{equation}\label{Mpro}
P_\lambda^{(l)}\ := \
-\D\frac{1}{2\pi\ri}\Int_\gamma dz\,\,[M_1(z,\Gamma_l)]^{-1}
\end{equation}
where $\gamma$ stands for an arbitrary rectifiable closed
contour going around $\lambda$ in the positive direction
in a sufficiently close neighbourhood such that
$\gamma\cap\Gamma_l=\emptyset$
and no points of the spectrum of $M_1(\cdot,\Gamma_l)$, except the
eigenvalue $\lambda$, lie inside $\gamma$. Then the following
relations hold:
\begin{equation}\label{MresiduePP}
P_\lambda^{(l)}\ =\ 
\sP_\lambda^{(l)}\,\,(I_1+\Omega^{(l)})^{-1}
\ =\ (I_1+\Omega^{(l)})^{-1}\,\,\sP_\lambda^{(-l)*}\,.
\end{equation}
\end{thm}
\begin{proof} The proof is carried out in the same way as the proof
of the relation~(\ref{MOmega}),
only the path of integration is changed.
\end{proof}
%%%%%%%%%%%%%%%%%%%%%%%%%%%%%%%%%%%%%%%%%%%%%%%%%%%%%%%%%%%%%%%%
\section{An example}\label{Example}
%%%%%%%%%%%%%%%%%%%%%%%%%%%%%%%%%%%%%%%%%%%%%%%%%%%%%%%%%%%%%%%%%%%
Let $\cH_0=\cH_1=L_2(\R)$ and $A_0=\sD^2+\lambda_0 I_0$ where
$\sD=\ri\D\frac{d}{dx}$ and $\lambda_0$ is some positive number.
It is assumed that the domain $\cD(A_0)$ is the Sobolev space
$W_2^2(\R)$. The spectrum of $A_0$ is absolutely continuous and
fills the semiaxis $\Delta _0=[\lambda_0,+\infty)$. By the
operator $A_1$ we understand
the multiplication by a bounded real-valued function
$a_1$, $A_1 f_1=a_1f_1$,
$f_1\in\cH_1$. The operator $T_{01}$ reads as
$$
   T_{01}=(\sD^2+\lambda_0 I_0)^{1/2}B
$$
where $B$ is the multiplication by a bounded
real-valued function $b\in W_2^1(\R)$, $Bf=bf$,
$f{\in}L_2(\R)$. Moreover, we assume that the function $b$ is
decreasing at infinity at least exponentially, so
that for any $x\in\R$ the estimate
\begin{equation}\label{EstimB}
   |b(x)|\ \le \ c\,\exp(-\alpha|x|)
\end{equation}
holds with some $c>0$ and $\alpha>0$.  Finally, we assume that
the range of the function
$$
    \ta_1(x)\ =\ a_1(x)-b^2(x)
$$
is embedded in $[\lambda _0+\widetilde{c},\infty)$ with 
some $\widetilde{c}>0$. The
operator $\tA_1$ is the multiplication by the function
$\ta_1 $.

It is easy to check that the value $E^0(\mu)$ of the spectral function
of the operator $A_0=\sD^2+\lambda_0 I_0$ is represented by 
the integral operator
whose kernel reads
$$
  E^0(\mu;x,x')\ =\ \left\{\begin{array}{cl}
0 & \mbox{ if } \mu <\lambda _0, \\[2mm]
\D\frac{1}{\sqrt{2\pi}}
\int_{\lambda_0}^\mu d\nu\,\,
\D\frac{\cos[(\nu-\lambda_0)^{1/2}(x-x')]}
{(\nu-\lambda_0)^{1/2}}
 & \mbox{ if } \mu\geq \lambda _0.
\end{array}\right.
$$
Thus the derivative $K'_B(\mu)$ is also an integral operator in $L_2(\R)$.
Its kernel $K'_B(\mu;x,x')$ is only nontrivial for $\mu>\lambda_0$
and, moreover, for these $\mu$
$$
K'_B(\mu;x,x')\ =\ \D\frac{1}{\sqrt{2\pi}}\,
\D\frac{\cos[(\mu-\lambda_0)^{1/2}(x-x')]}
{(\mu-\lambda_0)^{1/2}}\,\, b(x)\,b(x').
$$
Obviously, this kernel is degenerate for $\mu >\lambda _0$,
\begin{equation}\label{KDegenerate}
K'_B(\mu;x,x')\ =\ \D\frac{1}{2\, \sqrt{2\pi}\, (\mu -\lambda _0)^{1/2}}\,
[\tb_+(\mu,x)\tb_-(\mu,x')+\tb_-(\mu,x)\tb_+(\mu,x')]
\end{equation}
where
$\tb_\pm(\mu,x)={\rm e}^{\pm\ri\,(\mu-\lambda_0)^{1/2}x}\,b(x)$.
{}From the assumption (\ref{EstimB}) on $b$ we conclude that
in the domain $\pm\Img\sqrt{\mu-\lambda_0}<\alpha$, i.\,e., inside
the parabola
\begin{equation}\label{HolDom}
   \Real\mu>\lambda_0-\alpha^2+\D\frac{1}{4\alpha^2}(\Img\mu)^2,
\end{equation}
the functions $\tb_\pm(\mu,\cdot)$ are elements of $L_2(\R)$. The
function $K'_B(\mu)$ admits an analytic continuation onto this
domain (cut along the interval $ \lambda _0-\alpha^2<\mu\le\lambda_0$)
as a holomorphic
function with values in $\bB(\cH_1,\cH_1)$ and
the equation~(\ref{KDegenerate}) implies that
$$
   \|K'_B(\mu)\|\ \leq\ \D\frac{1}{\sqrt{2\pi}}\,
   \frac{1}{|(\mu-\lambda_0)^{1/2}|}\,
   \| \tb_+(\mu,\cdot)\|\,\| \tb_-(\mu,\cdot)\|.
$$
Obviously, for real $\mu$ we have $\|
\tb_\pm(\mu,\cdot)\|=\|b\|$.  Since $\tA_1$ is bounded, one can
always choose a $K_B$-bounded contour $\Gamma$ lying in the
domain~(\ref{HolDom}). Indeed, for the $K_B$-boundedness of
the contour $\Gamma$  it is sufficient to have its infinite part
presented by an appropriate semi-infinite real interval. Thus,
if the function $b$ is sufficiently small in the sense that the
conditions~(\ref{Best}) hold, one can apply 
all the statements of the Section
\ref{SmainEq} and \ref{SecFactor}
to the corresponding transfer function $M_1(z,\Gamma)$.

\end{document}